\documentclass{amsart}
\usepackage{amsfonts,amsmath, amssymb, amsthm}
\usepackage{mathtools}
\usepackage{graphicx}
\usepackage{geometry}
\usepackage[all]{xy}
\usepackage{tikz-cd}
\usepackage{tikz}
\usetikzlibrary{patterns}
\usepackage{bookmark}
\usepackage{hyperref}

\theoremstyle{plain}
\newtheorem{thm}{Theorem}[section]

\newtheorem{prop}[thm]{Proposition}
\newtheorem{cor}[thm]{Corollary}

\theoremstyle{definition}
\newtheorem{defi}[thm]{Definition}

\theoremstyle{remark}
\newtheorem{rem}[thm]{Remark}

\providecommand{\Cb}{\mathbb{C}}
\providecommand{\Rb}{\mathbb{R}}
\providecommand{\Sb}{\mathbb{S}}
\newcommand{\Ac}{\mathcal{A}}
\newcommand{\Lc}{\mathcal{L}}
\newcommand{\Pc}{\mathcal{P}}
\newcommand{\Int}{\mathrm{Int}}
\DeclareMathOperator{\Conf}{Conf}

\usepackage{caption}
\allowdisplaybreaks[2]
\newcommand{\bigprec}{\mathop{\mathrel{\raisebox{-0.1ex}{\scalebox{1.5}{$\prec$}}}}}

\title{A lexicographic section of the braid arrangement and the modified Artin presentation}
\author{So Yamagata}
\address{Department of Applied Mathematics, Fukuoka University}
\email{so.yamagata@fukuoka-u.ac.jp}
\subjclass[2020]{20F36, 52C35, 20F05, 14F35}
\keywords{pure braid group, braid arrangement, hyperplane arrangement, lexicographic section, braid monodromy, Manin--Schechtman arrangement}
\date{}

\begin{document}
\begin{abstract}
  We study a specific line arrangement obtained from a generic $2$--section of the braid arrangement, and compute the fundamental group of its complement via braid monodromy.
  We show that the resulting presentation of the fundamental group coincides, under the identification of generators, with the modified Artin presentation introduced by Margalit and McCammond.
  Moreover, we extend the construction to the Manin--Schechtman arrangements $MS(n, k)$, which are higher analogues of the braid arrangement.
  Focusing on the case $k = 2$, we obtain an explicit presentation of $\pi_1(\Cb^n \setminus MS(n, 2))$.
\end{abstract}

\maketitle

\section{Introduction}\label{sec:intro}
In 1891 Hurwitz~\cite{Hur91} introduced the braid group as the fundamental group of the configuration space of $n$ points in the complex plane.
Artin~\cite{Art25, Art47} subsequently established the algebraic and geometric foundations of braid groups and gave the classical presentations of both the braid group and the pure braid group.
While the braid group admits a relatively compact and tractable presentation, the situation for the pure braid group is markedly different.
As Artin famously observed:
\begin{quote}
  \textit{Although it has been proved that every braid can be deformed into a similar normal form the writer is convinced that any attempt to carry this out on a living person would only lead to violent protests and discrimination against mathematics. He would therefore discourage such an experiment.}
\end{quote}
\rightline{Emil Artin~\cite{Art47}}
This remark illustrates the intrinsic complexity of explicit calculations in the pure braid group and motivates the search for descriptions that are both explicit and geometrically meaningful.

Margalit and McCammond~\cite{MM09} made an attempt at such presentations.
They gave several positive finite presentations of the (pure) braid group of easy and simple form.
Among such presentations, we focus on what they called the modified Artin presentation.
Roughly speaking, the modified Artin presentation is obtained by considering the pure braid group as the mapping class group of a convexly punctured disc, with generators given by swings of the punctures and relations classified into three geometrically defined types.

From the viewpoint of hyperplane arrangements, the pure braid group arises as the fundamental group of the complement of the braid arrangement
\begin{equation*}
  Br(n) = \{ H_{ij} \mid 1 \leq i < j \leq n \}, \quad  H_{ij} = \{(x_1, x_2, \dots, x_n) \in \Cb^n \mid x_i = x_j \}.
\end{equation*}
Thus one has an identification
\begin{equation*}
  PB_n \cong \pi_1(\Cb^n \setminus Br(n)).
\end{equation*}

The topology of complements of hyperplane arrangements has been extensively studied from both combinatorial and geometric viewpoints (see, for instance~\cite{OT92, FR87, FR00, Yos24}).
In this paper, we apply the braid monodromy framework to carefully chosen generic $2$--sections of the braid arrangements.

A key input for applying this method to hyperplane arrangements is the Zariski--Lefschetz hyperplane section theorem (\cite{HuT73}), which implies that the fundamental group of the complexified real line arrangement is isomorphic to that of a generic $2$--dimensional section.
Consequently, one may replace the hyperplane arrangement by a suitable complexified real line arrangement arising from a generic $2$--section and compute the desired fundamental group using braid monodromy.
Such generic sections, however, are far from unique, and different choices may lead to presentations that vary significantly in their transparency and computational usefulness.
In order to make braid monodromy computations genuinely explicit, it is therefore necessary to choose a section for which the order in which intersection points appear under a horizontal sweep is controlled in a systematic way.

In this paper we introduce a lexicographic line arrangement $\Lc_n$ which is obtained as the generic $2$--section of the braid arrangement $Br(n)$.
Its defining feature is that the lines $l_{ij}$ and their triple intersection points $P_{ijk}$ occur in an order compatible with the lexicographic order of their indices.
This additional structure makes the wiring diagram of $\Lc_n$ computable by a uniform combinatorial rule and allows the Moishezon--Teicher algorithm to be carried out explicitly.
As a result, one obtains a concrete presentation of the pure braid group directly from the geometry of the arrangement.

One of the main results of the paper is that the presentation obtained from this braid monodromy computation coincides, under the identification of generators, with the modified Artin presentation provided in \cite{MM09}.
Although both the braid monodromy presentation and the modified Artin presentation are well known individually, the relationship between them does not appear to have been made explicit in existing work.
The lexicographic section thus provides a geometric explanation for the modified Artin presentation and shows that it arises naturally from a natural arrangement--theoretic construction.

The second result of the paper concerns a higher-dimensional analogue of the braid arrangement.
Manin and Schechtman~\cite{MS89} introduced a family of arrangements $MS(n, k)$ associated with a generic central arrangement of $n$ hyperplanes in $\Cb^k$.
The ambient parameter space may be identified with $\Cb^n$, and the hyperplanes of $MS(n, k)$ record degenerations of $(k + 1)$-subarrangements, namely those parameter values for which the corresponding translated hyperplanes fail to be in general position.
When $k = 1$ this construction recovers the braid arrangement, so that $MS(n, 1) = Br(n)$.
For $k \geq 2$ the arrangement exhibits higher-multiplicity intersection points governed by the combinatorics of $(k + 2)$-subsets, as described in~\cite{MS89}.
These arrangements provide a natural testing ground for extending braid monodromy techniques beyond the classical braid arrangement.

We generalize the construction of the lexicographic line arrangement to the Manin--Schechtman arrangements $MS(n, k)$.
More precisely, we define $(n, k)$--lexicographic line arrangements as generic $2$--sections of $MS(n, k)$ equipped with a compatible lexicographic ordering of their lines $l_{I}, I \subset [n], |I| = k + 1 $, and their multiple intersection points (of multiplicity $k + 2$) $P_{L}, L \subset [n], |L| = k + 2$.
In this paper we carry out the braid monodromy computation explicitly only in the case $k=2$.
As an application, we derive an explicit presentation of $\pi_1(\Cb^n \setminus MS(n, 2))$, thereby treating the full family of Manin--Schechtman arrangements in the case $k = 2$.
The case $k = 2$ already exhibits the new phenomena arising from higher--multiplicity intersection points beyond the braid arrangement.
More precisely, we show that the fundamental group of the complement of $MS(n, 2)$ admits a presentation that can be regarded as a natural generalization of the modified Artin presentation for the pure braid group.
Note that the case $k = 2$ is the minimal instance in which higher--multiplicity intersection points occur, and no essentially new local configurations appear for $k > 2$.

The paper is organized as follows.
Section~\ref{sec:modified_artin} reviews the modified Artin presentation of the pure braid group following~\cite{MM09}.
Section~\ref{sec:braid_monod} recalls the braid monodromy technique for complexified real line arrangements.
In Section~\ref{sec:braid_arr} we introduce the lexicographic line arrangement associated with the braid arrangement, compute its braid monodromy, and show that the resulting presentation agrees with the modified Artin presentation.
Finally, Section~\ref{sec:MS_pi1} discusses Manin--Schechtman arrangements, introduces their $(n, k)$-lexicographic line arrangements, and carries out the computation for $\pi_1(\Cb^n \setminus MS(n, 2))$.
\section{Modified Artin presentation for the pure braid group}\label{sec:modified_artin}
In this section, following \cite{MM09}, we recall the pure braid group as the fundamental group of a configuration space and state the modified Artin presentation.
We use convex swings as generators, and the defining relations are classified into three types according to configurations of punctures.
\begin{defi}
  A \textbf{convexly punctured disc} is a closed disc $D \subset \Rb^2$ equipped with punctures
  \begin{equation*}
    K = \{ a_1, a_2, \dots, a_n \} \subset \Int(D),
  \end{equation*}
  arranged so that every puncture lies on the boundary of the convex hull of $K$.
\end{defi}
Note that in this definition, the punctures are arranged so that they form the vertices of a convex $n$--gon.
This determines a cyclic order on $K$ along the polygonal boundary.
We use this cyclic order throughout; the particular geometric embedding is irrelevant up to isotopy.
Throughout the paper, we simply write $i$ for the puncture $a_i$ whenever no confusion arises.
In particular, we fix the labeling so that the punctures are sorted in ascending order in a counterclockwise manner.

\begin{defi}
  Let $A \subset K$.
  The \textbf{convexly punctured subdisc} $D_A$ is a sufficiently small neighborhood of the convex hull of the punctures in $A$, chosen so that $D_A \cap K = A$.
  Since the punctures are in convex position, such a neighborhood exists and is unique up to isotopy.
\end{defi}

\begin{defi}
  For a convexly punctured disc $D_K$, the \textbf{pure braid group} $PB_K$ is defined as the fundamental group of the configuration space $\Conf(D_K, |K|)$.
\end{defi}

\begin{defi}
  Let $A \subset K$, and let $P_A$ be the convex hull of $A$.
  A \textbf{convex rotation} $R_A$ is the rigid motion that moves each puncture in $A$ along the boundary of $P_A$ through one edge, preserving the cyclic order, while fixing all punctures outside $A$.
  The pure braid group condition requires each puncture to return to itself.
  Thus the smallest positive power of $R_A$ lying in $PB_K$ is
  \begin{equation*}
    S_A := (R_A)^{|A|},
  \end{equation*}
  called the \textbf{convex swing} associated with $A$.
\end{defi}
The following combinatorial conditions are used to describe the relations of the pure braid group.
\begin{defi}
  Two subsets $A,B\subset K$ are said to be \textbf{non-crossing} if the convex hulls of $A$ and $B$ are disjoint in their interiors; equivalently, the subdiscs $D_A$ and $D_B$ are disjoint.
  An ordered partition $(A_1, A_2, \dots, A_k)$ of $K$ is \textbf{admissible} if, with respect to the cyclic order on $K$, the points belonging to each $A_i$ appear as a single consecutive block, and these blocks occur in the order $A_1, A_2, \dots, A_k$.
\end{defi}

Having introduced convex swings as the basic geometric motions in a convexly punctured disc, we now turn to their role in generating the pure braid group.
In \cite{MM09}, Margalit and McCammond observed that every interaction among convex swings is governed by exactly three geometric configurations: non-crossing pairs, crossing pairs, and admissible triples of punctures.
Each configuration determines a corresponding commutation relation among the swings, and together these three types of relations constitute a complete and minimal defining set for the pure braid group.

\begin{thm}[Theorem~2.3 in \cite{MM09}]\label{thm:MM}
  Let $D_K$ be a convexly punctured disc.
  Then $PB_K$ is generated by the convex swings $S_{ij}$, and all defining relations among them fall into the following three types (see Figure~\ref{fig:puc_conf}):
  \begin{enumerate}
    \item[(R-I)]
          If the pairs $\{ i, j \}$ and $\{ r, s \}$ are non-crossing, then
          \begin{equation*}
            [S_{ij}, S_{rs}] = 1.
          \end{equation*}
    \item[(R-II)]
          If the pairs $\{ i, j \}$ and $\{ r, s \}$ cross in the counterclockwise order with $i, s, j, r$ in this order, then the corresponding swings satisfy
          \begin{equation*}
            [S_{ij}, S_{is} S_{rs} S_{is}^{-1}] = 1.
          \end{equation*}
    \item[(R-III)]
          If $(i, j, k)$ is an admissible triple (i.e., $i,j,k$ occur as a single consecutive block in this order in the cyclic ordering of $K$), then
          \begin{equation*}
            [S_{ij}, S_{ik} S_{jk}] = 1, \quad [S_{ij} S_{ik}, S_{jk}] = 1.
          \end{equation*}
  \end{enumerate}
\end{thm}
\begin{figure}[h]
  \centering
  \begin{minipage}[b]{0.29\linewidth}
    \centering
    \begin{tikzpicture}
      \draw (0,0) circle (1.7);
      \fill (-1,0) circle (3pt);
      \fill (0,-1) circle (3pt);
      \fill (1,0) circle (3pt);
      \fill (0,1) circle (3pt);
      \node at (-1.2, 0.2) {$i$};
      \node at (0.2, -1.2) {$j$};
      \node at (1.2, -0.2) {$r$};
      \node at (-0.2, 1.2) {$s$};
      \draw (-1,0) -- (0,-1);
      \draw (1,0) -- (0,1);
    \end{tikzpicture}
  \end{minipage} \
  \begin{minipage}[b]{0.29\linewidth}
    \centering
    \begin{tikzpicture}
      \draw (0,0) circle (1.7);
      \fill (-1,0) circle (3pt);
      \fill (0,-1) circle (3pt);
      \fill (1,0) circle (3pt);
      \fill (0,1) circle (3pt);
      \node at (0, 1.3) {$j$};
      \node at (0, -1.3) {$i$};
      \node at (-1.3, 0) {$r$};
      \node at (1.3, 0) {$s$};
      \draw (-1,0) -- (1,0);
      \draw (0,-1) -- (0,1);
    \end{tikzpicture}
  \end{minipage} \
  \begin{minipage}[b]{0.29\linewidth}
    \centering
    \begin{tikzpicture}
      \draw (0,0) circle (1.7);
      \fill (1,0) circle (3pt);
      \fill (0,1) circle (3pt);
      \fill (0,-1) circle (3pt);
      \node at (0,1.3) {$i$};
      \node at (0,-1.3) {$j$};
      \node at (1.3, 0) {$k$};
      \draw (1,0) -- (0,1) -- (0,-1) -- cycle;
    \end{tikzpicture}
  \end{minipage}
  \caption{The configurations of punctures corresponding to (R-I)--(R-III), from left to right.}
  \label{fig:puc_conf}
\end{figure}
\section{Braid monodromy technique}\label{sec:braid_monod}
In this section, following \cite{GTV02}, we briefly recall the notions of a line arrangement, its associated wiring diagram, and the algorithm used to compute the fundamental group of the complement of a complexified line arrangement.
\subsection{Line arrangement and its wiring diagram}
A \textbf{line arrangement} in $\Cb^2$ is a finite collection of complex lines.
The arrangement is called \textbf{real} if each line is defined by an equation with real coefficients.
Via a suitable affine identification of a real plane $\Pc \cong \Rb^2$ inside $\Cb^2$, a real arrangement gives rise to a line arrangement in the usual $(x, y)$-plane, which will serve as the input for the braid monodromy algorithm.

To a real arrangement $\Lc = \{ l_1, l_2, \dots, l_n \}$ one associates a \textbf{wiring diagram} \cite{Goo80}, which encodes the combinatorics of its intersection points.
After an orientation-preserving affine transformation, we may assume that no line of $\Lc$ is vertical and that no two intersection points have the same $x$-coordinate.
Under these assumptions, the wiring diagram associated with $\Lc$ is a collection of $n$ wires in $\Rb^2$ (each wire being a union of segments and rays, homeomorphic to $\Rb$).

For each $t \in \Rb$, let $l_t = \{ x = t \}$ be the vertical line at position $t$.
When $t \gg 0$, the intersections $l_t \cap l_i$ occur in a fixed vertical order with respect to the $y$-axis, which we take as the initial ordering of the wires.
As $t$ decreases, this order changes only when $l_t$ passes the $x$-coordinate of an intersection point of multiplicity $k \geq 2$.
At such a moment, the $k$ wires meeting at that point undergo a local move:
for $k = 2$ the two adjacent wires swap, while for $k \geq 3$ the entire block of $k$ adjacent wires reverses its vertical order (equivalently, is turned upside down).
The resulting pseudolines are precisely the wires of the wiring diagram.

Each elementary switch is recorded by a \textbf{Lefschetz pair}, namely the minimal and maximal indices of the two wires involved.
Conversely, starting from a list of Lefschetz pairs, one reconstructs the wiring diagram uniquely by reversing this procedure.

Figures~\ref{fig:exlin} and \ref{fig:exwire} show an arrangement of six lines and its wiring diagram.
This arrangement is an example of a lexicographic line arrangement of the braid arrangement $Br(4)$, to which we will refer later.
In the example, each line is labeled as $l_{ij}$ for $1 \leq i < j \leq 4$, and the intersection points are labeled by
\begin{equation*}
  P_{ij \mid kl} = l_{ij} \cap l_{kl},
  \qquad
  P_{ijk} = l_{ij} \cap l_{jk} \cap l_{ik} \quad (1 \leq i < j < k \leq 4).
\end{equation*}

\begin{figure}[h]
  \centering
  \begin{tikzpicture}[scale=1.5]
    \coordinate (0) at (1, 2);
    \coordinate (1) at (1.5, 2);
    \coordinate (2) at (2, 1);
    \coordinate (3) at (2, 0.625);
    \coordinate (4) at (1, -1);
    \coordinate (5) at (-0.6, -2.1);
    \coordinate (6) at (1.5, 0.75);
    \coordinate (7) at (0.5, 1);
    \coordinate (8) at (0, 0);
    \coordinate (9) at (-0.25, 0.25);
    \coordinate (10) at (-0.625, -1.5);
    \coordinate (11) at (-1, -0.5);
    \coordinate (12) at (-1.5, 1.5);
    \coordinate (13) at (-1.6, 1.9);
    \coordinate (14) at (-2, 2);
    \coordinate (15) at (-2, 1.625);
    \coordinate (16) at (-2, -1);
    \coordinate (17) at (-2, -1.5);
    \coordinate (18) at (-1, -2);
    \coordinate [label=$l_{12}$] (a) at (1.1,2);
    \coordinate [label=$l_{13}$] (b) at (1.6, 1.9);
    \coordinate [label=$l_{14}$] (c) at (2.2, 0.9);
    \coordinate [label=$l_{23}$] (d) at (2.2, 0.4);
    \coordinate [label=$l_{24}$] (e) at (1.2, -1.3);
    \coordinate [label=$l_{34}$] (f) at (-0.5, -2.5);
    \coordinate [label={$P_{14 \mid 23}$}] (A) at (1.4, 0.8);
    \coordinate [label={$P_{123}$}] (B) at (0.3, 1.0);
    \coordinate [label={$P_{124}$}] (C) at (0.3, -0.2);
    \coordinate [label={$P_{13 \mid 24}$}] (D) at (-0.7, 0.1);
    \coordinate [label={$P_{12 \mid 34}$}] (E) at (-0.4, -1.7);
    \coordinate [label={$P_{134}$}] (F) at (-1.3, -0.6);
    \coordinate [label={$P_{234}$}] (G) at (-1.2, 1.4);
    \coordinate (H) at (3/2, 0.75);
    \coordinate (I) at (-1/4, 0.25);
    \coordinate (J) at (-3/4, -1.5);
    \begin{scope}
      \draw (0) -- (18);
      \draw (1) -- (17);
      \draw (2) -- (16);
      \draw (3) -- (15);
      \draw (4) -- (14);
      \draw (5) -- (13);
      \fill (7) circle [radius=0.05cm];
      \fill (8) circle [radius=0.05cm];
      \fill (11) circle [radius=0.05cm];
      \fill (12) circle [radius=0.05cm];
      \fill (H) circle [radius=0.05cm];
      \fill (I) circle [radius=0.05cm];
      \fill (J) circle [radius=0.05cm];
    \end{scope}
  \end{tikzpicture}\caption{The line arrangement of six lines with seven intersection points.}\label{fig:exlin}
\end{figure}

\begin{figure}[h]
  \centering
  \begin{tikzpicture}
    \node at (11.5,1) {$l_{34}$};
    \node at (11.5,2) {$l_{24}$};
    \node at (11.5,3) {$l_{23}$};
    \node at (11.5,4) {$l_{14}$};
    \node at (11.5,5) {$l_{13}$};
    \node at (11.5,6) {$l_{12}$};
    \node at (1,-0.3) {$P_{234}$};
    \node at (2.5,-0.3) {$P_{134}$};
    \node at (4,-0.3) {$P_{12 \mid 34}$};
    \node at (5.5,-0.3) {$P_{13 \mid 24}$};
    \node at (7,-0.3) {$P_{124}$};
    \node at (8.5,-0.3) {$P_{123}$};
    \node at (10,-0.3) {$P_{14 \mid 23}$};
    \draw (0,0) -- (11,0);
    \draw[dotted] (1,5) -- (1,0);
    \draw[dotted] (2.5,3) -- (2.5,0);
    \draw[dotted] (4,1.5) -- (4,0);
    \draw[dotted] (5.5,4.5) -- (5.5,0);
    \draw[dotted] (7,3) -- (7,0);
    \draw[dotted] (8.5,5) -- (8.5,0);
    \draw[dotted] (10,3.5) -- (10,0);
    \draw (0,1) -- (3.8,1) -- (4.2,2) -- (6.8,2) -- (7.2,4) -- (8.3,4) -- (8.7,6) -- (11,6);
    \draw (0,2) -- (2.3,2) -- (2.7,4) -- (5.3,4) -- (5.7,5) -- (11,5);
    \draw (0,3) -- (9.8,3) -- (10.2,4) -- (11,4);
    \draw (0,4) -- (0.8,4) -- (1.2,6) -- (8.3,6) -- (8.7,4) -- (9.8,4) -- (10.2,3) -- (11,3);
    \draw (0,5) -- (5.3,5) -- (5.7,4) -- (6.8,4) -- (7.2,2) -- (11,2);
    \draw (0,6) -- (0.8,6) -- (1.2,4) -- (2.3,4) -- (2.7,2) -- (3.8,2) -- (4.2,1) -- (11,1);
  \end{tikzpicture}\caption{The wiring diagram of the line arrangement in Figure~\ref{fig:exlin}.}\label{fig:exwire}
\end{figure}

\subsection{Computation of the fundamental group}\label{subsec:comp}
We recall the method for computing the fundamental group $\pi_1(\Cb^2 \setminus \Lc)$ of the complement of a complexified real line arrangement $\Lc$ from the data of its wiring diagram; see, for instance, \cite{GTV02, MT88} for detailed proofs.

Let $l = |\Lc|$ denote the number of lines in the arrangement.
Let $D = \{ z \in \Cb \mid |z - \tfrac{l + 1}{2}| \leq \tfrac{l + 1}{2} \}$ be a disc with punctures $K = \{ 1, 2, \dots, l \}$, and fix a base point $u_0 \in D$ below the real line.
Recall that the braid group $B_l$ arises as the mapping class group of $(D, K)$, i.e., diffeomorphisms of $D$ fixing the boundary pointwise and permuting $K$.
The pure braid group corresponds to the subgroup, where each puncture returns to itself.
The fundamental group $\pi_1(D \setminus K, u_0)$ is freely generated by $\{ g_1, g_2, \dots, g_l \}$, where $g_i$ is a simple loop based at $u_0$ encircling the $i$-th puncture once in the positive direction.

Each intersection point of the arrangement determines a \textbf{skeleton}, namely an embedded path in $D\setminus K$ that encodes the local monodromy around that point.
For the Lefschetz pair $[a_i, b_i]$, the \textbf{initial skeleton} is chosen to be a simple arc in $D \setminus K$ joining the consecutive punctures from $a_i$ to $b_i$.

If $[a_1, b_1], [a_2, b_2], \dots, [a_{i - 1}, b_{i - 1}]$ are the preceding pairs, listed from right to left, then the skeleton $s_i$ associated to the $i$-th intersection point is obtained by applying to the initial arc $s_i^{\mathrm{init}}$ the clockwise $\pi$-rotations corresponding to these pairs.
More precisely, for a Lefschetz pair $[a_j, b_j]$ let $\rho_{[a_j, b_j]}$ denote the clockwise half twist supported on the subdisc whose boundary contains exactly the punctures with labels between $a_j$ and $b_j$, fixing all other punctures pointwise.
Then, we obtain the skeleton
\begin{equation*}
  s_i =
  \rho_{[a_{1}, b_{1}]} \;
  \rho_{[a_{2}, b_{2}]} \dots
  \rho_{[a_{i - 1}, b_{i - 1}]} (s_i^{\mathrm{init}}).
\end{equation*}

For example, Figure~\ref{fig:skeleton_s3} illustrates the action of Lefschetz pairs on the initial skeleton $s_3^{\mathrm{init}}$ corresponding to the third Lefschetz pair $[3, 5]$ from the rightmost in Figure~\ref{fig:exwire}.
Each point $i$ in the figure corresponds to the point $C_{ij}=l_{ij}\cap l_t$ in Figure~\ref{fig:exlin}, where $l_t$ is the sectional line at position $t$.

\begin{figure}[h]
  \centering
  \begin{minipage}[b]{0.18\textwidth}
    \centering
    \begin{tikzpicture}
      \fill (0,0) circle (3pt);
      \fill (0,1) circle (3pt);
      \fill (0,2) circle (3pt);
      \fill (0,3) circle (3pt);
      \fill (0,4) circle (3pt);
      \fill (0,5) circle (3pt);
      \node at (0.4,0) {6};
      \node at (0.4,1) {5};
      \node at (0.4,2) {4};
      \node at (0.4,3) {3};
      \node at (0.4,4) {2};
      \node at (0.4,5) {1};
      \draw (0,1) -- (0,2) -- (0,3);
    \end{tikzpicture}
  \end{minipage}
  \begin{minipage}[b]{0.18\textwidth}
    \centering
    \raisebox{1cm}{
      \begin{tikzpicture}
        \fill (0.5,0.5) circle (3pt);
        \fill (0.5,1) circle (3pt);
        \fill (0.5,1.5) circle (3pt);
        \draw[fill=white] (0.5,2) circle (3pt);
        \fill (0.5,2.5) circle (3pt);
        \draw[fill=white] (0.5,3) circle (3pt);
        \fill (0,2.7) circle (3pt);
        \fill (1,2.3) circle (3pt);
        \draw[->] (-0.5,0) -- (1.5,0);
        \draw (0.5,1) -- (0.5,1.5);
        \draw[->] (0.4,2.1) .. controls (0,2.6) .. (0.4,2.9);
        \draw[->] (0.6,2.8) .. controls (1,2.4) .. (0.6,2.1);
        \draw (0.5,1.5) .. controls (0,2) .. (0,2.7);
      \end{tikzpicture}
    }
  \end{minipage}
  \begin{minipage}[b]{0.2\textwidth}
    \centering
    \begin{tikzpicture}
      \fill (0,0) circle (3pt);
      \fill (0,1) circle (3pt);
      \fill (0,2) circle (3pt);
      \fill (0,3) circle (3pt);
      \fill (0,4) circle (3pt);
      \fill (0,5) circle (3pt);
      \node at (0.4,0) {6};
      \node at (0.4,1) {5};
      \node at (0.4,2) {4};
      \node at (0.4,3) {3};
      \node at (0.4,4) {2};
      \node at (0.4,5) {1};
      \draw (0,1) -- (0,2);
      \draw (0,2) .. controls (-0.4,3) and (-0.4,4) .. (0,5);
    \end{tikzpicture}
  \end{minipage}
  \begin{minipage}[b]{0.18\textwidth}
    \centering
    \raisebox{1cm}{
      \begin{tikzpicture}
        \fill (0.5,0.5) circle (3pt);
        \fill (0.5,1) circle (3pt);
        \draw[fill=white] (0.5,1.5) circle (3pt);
        \fill (0,2) circle (3pt);
        \fill (1,2) circle (3pt);
        \draw[fill=white] (0.5,2.5) circle (3pt);
        \fill (0.5,3) circle (3pt);
        \fill (0.5,3.5) circle (3pt);

        \draw[->] (-0.5,0) -- (1.5,0);
        \draw (0.5,1) .. controls (0.1,1.5) .. (0,2);
        \draw (0,2) .. controls (0.1,3) .. (0.5,3.5);
        \draw[->] (0.4,1.6) .. controls (0,2) .. (0.4,2.4);
        \draw[->] (0.6,2.4) .. controls (1,2) .. (0.6,1.6);
      \end{tikzpicture}
    }
  \end{minipage}
  \begin{minipage}[b]{0.18\textwidth}
    \centering
    \begin{tikzpicture}
      \fill (0,0) circle (3pt);
      \fill (0,1) circle (3pt);
      \fill (0,2) circle (3pt);
      \fill (0,3) circle (3pt);
      \fill (0,4) circle (3pt);
      \fill (0,5) circle (3pt);
      \node at (0.4,0) {6};
      \node at (0.4,1) {5};
      \node at (0.4,2) {4};
      \node at (0.4,3) {3};
      \node at (0.4,4) {2};
      \node at (0.4,5) {1};
      \draw (0,1) .. controls (-0.4,2) .. (0,3);
      \draw (0,3) .. controls (-0.4,4) .. (0,5);
    \end{tikzpicture}
  \end{minipage}\caption{Action of Lefschetz pairs on $s_3^{\mathrm{init}}$.}\label{fig:skeleton_s3}
\end{figure}

From the resulting skeleton one obtains a relation in $\pi_1(D \setminus K, u_0)$, which later descends to a relation in $\pi_1(\Cb^2 \setminus \Lc)$ via the Zariski--van Kampen theorem.
Following Moishezon--Teicher, we distinguish the cases according to the multiplicity of the intersection point.

Suppose first that the Lefschetz pair $[a_i, b_i]$ corresponds to a double point.
In this case the skeleton $s_i$ is an embedded arc connecting the punctures $a_i$ and $b_i$.
Choosing an interior point on $s_i$ and pulling it to the base point $u_0$, the arc decomposes into two paths from $u_0$ to $a_i$ and $b_i$, respectively.
Each path determines a loop in $\pi_1(D\setminus K,u_0)$ by going once positively around the corresponding puncture and returning along the same path.

The resulting relation is a commutator
\begin{equation*}
  [A_{a_i}, B_{b_i}] = 1,
\end{equation*}
where $A_{a_i}$ and $B_{b_i}$ denote the two loops associated with the endpoints of the skeleton.
While $A_{a_i}$ is conjugate to the standard generator $g_{a_i}$, the loop $B_{b_i}$ is in general a conjugate of $g_{b_i}$ by a word determined by the combinatorics of the skeleton.
More precisely, as the path from $u_0$ to $b_i$ crosses other punctures, the corresponding loop is conjugated by the generators associated with those punctures.
Thus one may write
\begin{equation*}
  B_{b_i} = w g_{b_i} w^{-1},
\end{equation*}
where the word $w$ records the relative position of the skeleton with respect to the intervening punctures.
Consequently, the relation obtained from a double point takes the form
\begin{equation*}
  [g_{a_i}, w g_{b_i} w^{-1}] = 1,
\end{equation*}
and this relation depends only on the isotopy class of the skeleton.

More generally, suppose that $[a_i, b_i]$ corresponds to a multiple intersection point of multiplicity $k$.
In this case the skeleton produces $k$ loops $A_1, A_2, \dots, A_k$ in $\pi_1(D\setminus K, u_0)$, each associated with one of the punctures involved.
As in the double point case, each $A_j$ is conjugate to the standard generator corresponding to its endpoint, with the conjugating word determined by how the skeleton passes the remaining punctures.

The braid monodromy relation is then expressed as a cyclic equality
\begin{equation*}
  A_1 A_2 \cdots A_k = A_2 A_3 \cdots A_k A_1 = \cdots = A_k A_1 \cdots A_{k-1},
\end{equation*}
which is equivalently written as a family of commutator relations
\begin{equation*}
  [A_1, A_2 \cdots A_k] = [A_1A_2, A_3 \cdots A_k] = \cdots = [A_1 \cdots A_{k-1}, A_k] = 1.
\end{equation*}
Each of these relations can be rewritten explicitly in terms of the free generators $\{g_1,\dots,g_l\}$ by substituting the corresponding conjugated expressions.
In particular, the precise form of the relations is completely determined by the combinatorial data of the skeleton, namely by the order in which it passes the punctures.
\section{Braid arrangement and braid monodromy}\label{sec:braid_arr}
In this section, using the algorithm recalled in Section \ref{sec:braid_monod}, we derive an explicit presentation of the fundamental group of the complement of the braid arrangement.
It is well known that this group is the pure braid group.
By the Zariski--Lefschetz hyperplane section theorem, it suffices to consider a generic $2$-section $\Pc$ of the arrangement, namely the associated complexified line arrangement.
\subsection{A lexicographic line arrangement}\label{subsec:lex}
We begin by recalling the definition of the braid arrangement.
\begin{defi}
  For $1 \leq i < j \leq n$, let
  \begin{equation*}
    H_{ij} = \{ (x_1, x_2, \dots, x_n) \in \Cb^n \mid x_i = x_j \}
  \end{equation*}
  be the hyperplane in $\Cb^n$.
  The \textbf{braid arrangement} $Br(n)$ is an arrangement of hyperplanes defined by
  \begin{equation*}
    Br(n) = \{ H_{ij} \}_{1 \leq i < j \leq n}.
  \end{equation*}
\end{defi}
A generic $2$-section of $Br(n)$ yields a complexified line arrangement in $\Cb^2$.
In this subsection, we construct a particular such section, which we call the lexicographic line arrangement and denote by $\Lc_n$.

Let $l_{ij} = H_{ij} \cap \Pc$ be the complex line obtained by restricting $H_{ij}$ to the generic complex $2$--plane $\Pc \cong \Cb^2$.
Thus $\{l_{ij}\}$ forms a complexified line arrangement in $\Pc$.
For visualization and for the purpose of applying the braid monodromy technique, we choose the section $\Pc$ so that the resulting arrangement is defined over $\Rb^2$, and identify $\Pc$ with $\Rb^2$ equipped with its standard real structure.
In this way, the family $\{ l_{ij} \}$ may be viewed as a real line arrangement in a real plane $\Rb^2$, which we denote by $\Lc_n$.

Since the braid arrangement has only two types of codimension--$2$ intersections,
\begin{equation*}
  H_{ij} \cap H_{kl} \quad \text{and} \quad H_{ij}\cap H_{ik}\cap H_{jk},
\end{equation*}
any generic $2$--section yields a line arrangement having only double points $l_{ij}\cap l_{kl}$ and triple points $l_{ij} \cap l_{ik} \cap l_{jk}$.

For a point $P \in \Rb^2$ we denote its $x$- and $y$-coordinates by $x(P)$ and $y(P)$.

\begin{defi}
  The \textbf{lexicographic line arrangement} $\Lc_n$ is a particular real line arrangement consisting of lines $l_{ij}$, $1 \leq i < j \leq n$, satisfying the following properties.

  Fix a circle $C$ that contains all intersection points of the lines $l_{ij}$ in its interior.
  For each line $l_{ij}$, let $C_{ij}$ denote the intersection point $C\cap l_{ij}$ lying on the ray of $l_{ij}$ in the direction of increasing $x$.
  Then,
  \begin{itemize}
    \item the points $C_{ij}$ appear in lexicographic order when going along $C$ clockwise starting from $C_{12}$;
    \item the only intersection points are of the following two types:
          \begin{itemize}
            \item double point $P_{ij|kl} = l_{ij} \cap l_{kl}$ for disjoint pairs $\{ i, j \}$ and $\{ k, l \}$;
            \item triple point $P_{ijk} = l_{ij} \cap l_{ik} \cap l_{jk}$ for $3$--subset $\{ i, j, k \} \subset [n]$;
          \end{itemize}
    \item the triple points are ordered by $x(P_{i'j'k'}) < x(P_{ijk})$ whenever $(i,j,k) < (i',j',k')$ as in the lexicographic order.
  \end{itemize}
\end{defi}

Figure~\ref{fig:exlin} shows an example of the lexicographic line arrangement $\Lc_4$.
The proof of the following proposition will be given later in Proposition \ref{prop:nk-line-exist} in a more general setting.
\begin{prop}
  The lexicographic line arrangement $\Lc_n$ exists.
\end{prop}
\subsection{Computation of braid monodromy and presentation for the pure braid group}\label{subsec:monodromy_tech}
In this subsection, we apply the braid monodromy technique to the lexicographic line arrangement $\Lc_n$ constructed in the previous subsection.
Thanks to the special ordering properties of $\Lc_n$, we can compute the braid monodromy explicitly and derive a presentation of the pure braid group $PB_n$.
We then show that this presentation agrees with the modified Artin presentation recalled in Section~\ref{sec:modified_artin}.

For each sectional line $l_t = \{ x = t \}$, where $-\infty < t < \infty$, let $Q_{ij}^t = l_{ij}\cap l_t$ denote the intersection point for $1 \leq i < j \leq n$.
As the line $l_t$ crosses an intersection point of the arrangement, the vertical order of the points $Q_{ij}^t$ changes locally.
To record this order, we write $(i,j) \prec (k,l)$ whenever $y(Q_{ij}^t) > y(Q_{kl}^t)$.

For $t$ sufficiently large, this order is lexicographic:
\begin{equation*}
  (1, 2) \prec (1, 3) \prec \dots \prec (n-1, n),
\end{equation*}
while for $t$ sufficiently negative, the order becomes colexicographic:
\begin{equation*}
  (n-1, n) \prec (n-2, n) \prec \dots \prec (1, 2),
\end{equation*}
We will describe how the vertical order changes as $l_t$ moves from $+\infty$ to $-\infty$;
before doing so, we introduce some notation.

For $s_1 \leq s_2 < t$ and $0 \leq i \leq s_2 - s_1 + 1$ we define an ordered sequence of pairs
\begin{equation*}
  ([s_1, s_2], t)^{(i)} \coloneqq (s_1, t) \prec (s_1 + 1, t) \prec \dots \prec
  \overset{\bar{i}}{(1, t)} \prec \dots \prec (s_2, t),
\end{equation*}
where the element $(s_1, t)$ is placed at the $i$-th position counted from the maximum element $(s_2, t)$.
If $i = 0$, we set
\begin{equation*}
  ([s_1, s_2], t)^{(0)} \coloneqq (s_1, t) \prec (s_1 + 1, t) \prec \dots \prec (s_2, t).
\end{equation*}
For convenience, we also allow concatenations of such blocks.
For fixed $i$ and fixed $s_1 \leq s_2$, and integers $t_1 \leq t_2$, we write
\begin{equation*}
  \bigprec_{\tau = t_1}^{t_2} ([s_1, s_2], \tau)^{(i)}
\end{equation*}
for the sequence
\begin{equation*}
  ([s_1, s_2], t_1)^{(i)} \prec ([s_1, s_2], t_1 + 1)^{(i)} \prec \dots \prec ([s_1, s_2], t_2)^{(i)}.
\end{equation*}

Finally, we introduce shorthand for lexicographic and colexicographic sequences.
For pairs $(s_1, t_1)$ and $(s_2, t_2)$ with $s_1 \leq s_2$ and $t_1 \leq t_2$, we set
\begin{align*}
  \mathrm{lex}\big((s_1, t_1),(s_2, t_2)\big)
   & := (s_1, t_1) \prec (s_1, t_1+1) \prec \dots \prec (s_2, t_2),
   &                                                                & \text{(lexicographic order)},   \\[4pt]
  \mathrm{colex}\big((s_1, t_1),(s_2, t_2)\big)
   & := (s_2, t_2) \prec (s_2, t_2-1) \prec \dots \prec (s_1, t_1),
   &                                                                & \text{(colexicographic order)}.
\end{align*}

The changes in the vertical order are generated by two elementary local operations:
\begin{itemize}
  \item $(i, j) \prec (k, l)\;\longleftrightarrow\;(k, l) \prec (i, j)$
        called a \textbf{single exchange (s.e.)}
  \item $(i, j) \prec (i, k) \prec (j, k)\;\longleftrightarrow\;(j, k) \prec (i, k) \prec (i, j)$  called a \textbf{flip}.
\end{itemize}
By combining these two operations, one obtains the entire evolution of the order as $l_t$ moves from $+\infty$ to $-\infty$, and hence the wiring diagram.

We now describe the order changes in detail.
First, starting from sufficiently large $t$ so that the order is lexicographic, we sweep the line $l_t$ in the negative $x$--direction and record the order changes up to the first triple point encountered, namely $P_{1,n-1,n}$.

(I) Order changes for $x(P_{123}) - \epsilon < t < + \infty$
\begin{align*}
                                            & (1, 2) \prec (1, 3) \prec \dots \prec (1, n - 1) \prec {\color{red}{(1, n) \prec (2, 3)}} \prec (2, 4) \prec \dots \prec (n-1, n) \\
  \xrightarrow{\text{s.e.'s with $(2, 3)$}} & {\color{red}{(1, 2) \prec (1, 3) \prec (2, 3)}} \prec (1, 4) \prec \dots \prec (1, n) \prec (2, 4) \prec \dots \prec (n-1, n)     \\
  \xrightarrow{\text{flip}}                 & (2, 3) \prec (1, 3) \prec (1, 2) \prec (1, 4) \prec \dots \prec (1, n) \prec (2, 4) \prec \dots \prec (n-1, n).
\end{align*}
Here, ``s.e.'s with $(2, 3)$'' means successive single exchanges of $(2, 3)$ with its adjacent elements in the current order.

(II) Order changes when $x(P_{1jk+1}) - \epsilon < t < x(P_{1jk}) - \epsilon$, where $2 \leq j < k \leq n - 1$.
\begin{align*}
                                            & \bigprec_{3 \leq t \leq j} ([2, t-1], t)^{(0)} \prec \bigprec_{j+1 \leq t \leq k} ([2, j], t)^{(1)} \prec {\color{red}{(1, j) \prec {\color{red}{(2, k+1)}}}} \prec ([3, j], k+1)^{(2)} \\
                                            & \prec \bigprec_{k+2 \leq t \leq n} ([2, j], t)^{(2)} \prec \mathrm{lex}((j+1, j+2), (n-1, n)) \prec \mathrm{colex}((1, 2), (1, j-1))                                                    \\
  \xrightarrow{\text{s.e.'s with $(1, j)$}} & \bigprec_{3 \leq t \leq j} ([2, t-1], t)^{(0)} \prec \bigprec_{j+1 \leq t \leq k} ([2, j], t)^{(1)}  \prec ([2, j-1], k+1)^{(0)}                                                        \\
                                            & \prec {\color{red}{(1, j) \prec (1, k+1) \prec (j, k+1)}} \prec \bigprec_{k+2 \leq t \leq n} ([2, j], t)^{(2)} \prec \mathrm{lex}((j+1, j+2), (n-1, n))                                 \\
                                            & \prec \mathrm{colex}((1, 2), (1, j-1))                                                                                                                                                  \\
  \xrightarrow{\text{flip}}                 & \bigprec_{3 \leq t \leq j} ([2, t-1], t)^{(0)} \prec \bigprec_{j+1 \leq t \leq k} ([2, j], t)^{(1)}  \prec ([2, j-1], k+1)^{(0)}                                                        \\
                                            & \prec (j, k+1) \prec (1, k+1) \prec (1, j) \prec \bigprec_{k+2 \leq t \leq n} ([2, j], t)^{(2)} \prec \mathrm{lex}((j+1, j+2), (n-1, n))                                                \\
                                            & \prec \mathrm{colex}((1, 2), (1, j-1)).
\end{align*}

(III) Order changes when $x(P_{1j+1j+2}) - \epsilon < t < x(P_{1jn}) - \epsilon$. \\
\begin{align*}
                                                & \bigprec_{3 \leq t \leq j} ([2, t-1], t)^{(0)} \prec \bigprec_{j+1 \leq t \leq n} ([2, j], t)^{(1)} \prec {\color{red}{(1, j) \prec \mathrm{lex}((j+1, j+2), (n-1, n))}} \\
                                                & \prec \mathrm{colex}((1, 2), (1, j-1))                                                                                                                                   \\
  \xrightarrow{\text{s.e.'s with $(1, j)$}}     & \bigprec_{3 \leq t \leq j} ([2, t-1], t)^{(0)} \prec \bigprec_{j+1 \leq t \leq n} ([2, j], t)^{(1)} \prec \mathrm{lex}((j+1, j+2), (n-1, n))                             \\
                                                & \prec \mathrm{colex}((1, 2), (1, j))                                                                                                                                     \\
  =                                             & \bigprec_{3 \leq t \leq j} ([2, t-1], t)^{(0)} \prec ([2, j], j+1)^{(0)} \prec \color{red}{(1, j+1) \prec ([2, j], j+2)^{(0)} \prec (1, j+2)}                            \\
                                                & \prec \bigprec_{j+3 \leq t \leq n} ([2, j], t)^{(1)} \prec \mathrm{lex}((j+1, j+2), (n-1, n)) \prec \mathrm{colex}((1, 2), (1, j))                                       \\
  \xrightarrow{\text{s.e.'s with $(1, j+1)$}}   & \bigprec_{3 \leq t \leq j+1} ([2, t-1], t)^{(0)} \prec ([2, j], j+2)^{(0)} \prec (1, j+1) \prec (1, j+2)                                                                 \\
                                                & \prec {\color{red}{\bigprec_{j+3 \leq t \leq n} ([2, j], t)^{(1)}}} {\color{red}{\prec (j+1, j+2)}} \prec \mathrm{lex}((j+1, j+3), (n-1, n))                             \\
                                                & \prec \mathrm{colex}((1, 2), (1, j))                                                                                                                                     \\
  \xrightarrow{\text{s.e.'s with $(j+1, j+2)$}} & \bigprec_{3 \leq t \leq j+1} ([2, t-1], t)^{(0)} \prec ([2, j], j+2)^{(0)} \prec {\color{red}{(1, j+1) \prec (1, j+2) \prec (j+1, j+2)}}                                 \\
                                                & \prec \bigprec_{j+3 \leq t \leq n} ([2, j], t)^{(1)} \prec \mathrm{lex}((j+1, j+3), (n-1, n)) \prec \mathrm{colex}((1, 2), (1, j))                                       \\
  \xrightarrow{\text{flip}}                     & \bigprec_{3 \leq t \leq j+1} ([2, t-1], t)^{(0)} \prec ([2, j], j+2)^{(0)} \prec (j+1, j+2) \prec (1, j+2) \prec (1, j+1)                                                \\
                                                & \prec \bigprec_{j+3 \leq t \leq n} ([2, j], t)^{(1)} \prec \mathrm{lex}((j+1, j+3), (n-1, n)) \prec \mathrm{colex}((1, 2), (1, j)).
\end{align*}

Next, we describe the order changes that occur as the sectional line $l_t$ moves further to the left of $x(P_{1n-1n})$, continuing the sweep in the negative $x$--direction.
In the computations below, we omit the final colexicographic block $\mathrm{colex}((1,2),(1,n))$ and record only the essential order changes.
To present the resulting expressions in a uniform form, we apply the index
relabeling
\begin{equation*}
  i \mapsto 1, \quad j \mapsto j - i + 1, \quad k \mapsto k - i + 1.
\end{equation*}
This change of variables is merely a relabeling of indices and does not affect the resulting relations.

(IV) Order changes when $x(P_{1jk+1}) - \epsilon < t < x(P_{1jk}) - \epsilon$.
\begin{align*}
   & \bigprec_{3 \leq t \leq j} ([2,t-1],t)^{(0)} \prec \bigprec_{j+1 \leq t \leq k-1} ([2,t-1],t)^{(t-j)} \prec
  ([2,j-1],k)^{(0)} \prec (j, k) \prec (1, k) \prec {\color{red}{(1, j)}}                                                                                                                                \\
   & {\color{red}{\prec ([j+1,k-1],k)^{(0)} \prec ([2, j-1], k+1)^{(0)}}} \prec (1, k+1) \prec (j, k+1) \prec ([j+1, k], k+1)^{(0)}                                                                      \\
   & \prec \bigprec_{k+2 \leq t \leq n} ([2,t-1],t)^{(t-j+1)} \prec \mathrm{colex}((1, 2), (1, j-1))                                                                                                     \\
   & \xrightarrow{\text{s.e.'s with $(1, j)$}} \bigprec_{3 \leq t \leq j} ([2,t-1],t)^{(0)} \prec \bigprec_{j+1 \leq t \leq k} ([2,t-1],t)^{(t-j)} \prec ([2,j-1],k+1)^{(0)} \prec {\color{red}{(1, j)}} \\
   & {\color{red}{\prec (1, k+1) \prec (j, k+1)}} \prec ([j+1,k],k+1)^{(0)} \prec \bigprec_{k+2 \leq t \leq n} ([2,t-1],t)^{(t-j+1)}                                                                     \\
   & \prec \mathrm{colex}((1, 2), (1, j-1))                                                                                                                                                              \\
   & \xrightarrow{\text{flip}} \bigprec_{3 \leq t \leq j} ([2,t-1],t)^{(0)} \prec \bigprec_{j+1 \leq t \leq k} ([2,t-1],t)^{(t-j)} \prec ([2,j-1],k+1)^{(0)} \prec (j, k+1) \prec (1, k+1)               \\
   & \prec (1, j) \prec ([j+1,k],k+1)^{(0)} \prec \bigprec_{k+2 \leq t \leq n} ([2,t-1],t)^{(t-j+1)} \prec \mathrm{colex}((1, 2), (1, j-1))                                                              \\
\end{align*}

(V) Order changes when $x(P_{1j+1j+2}) - \epsilon < t < x(P_{1jn}) - \epsilon$.

\begin{align*}
   & \bigprec_{3 \leq t \leq j} ([2,t-1],t)^{(0)} \prec \bigprec_{j+1 \leq t \leq n-1} ([2,t-1],t)^{(t-j)} \prec
  ([2,j-1],n)^{(0)} \prec (j, n) \prec (1, n) \prec {\color{red}{(1, j)}}                                                                                                                \\
   & {\color{red}{\prec ([j+1,n-1],n)^{(0)}}} \prec \mathrm{colex}((1, 2), (1, j-1))                                                                                                     \\
   & \xrightarrow{\text{s.e.'s with $(1, j)$}} \bigprec_{3 \leq t \leq j} ([2,t-1],t)^{(0)} \prec \bigprec_{j+1 \leq t \leq n} ([2, t-1],t)^{(t-j)} \prec \mathrm{colex}((1, 2), (1, j)) \\
   & = \bigprec_{3 \leq t \leq j} ([2,t-1],t)^{(0)} \prec ([2,j],j+1)^{(0)} \prec {\color{red}{(1, j+1) \prec ([2,j],j+2)^{(0)} }}                                                       \\
   & \prec (1, j+2) \prec (j+1, j+2) \prec \bigprec_{j+3 \leq t \leq n} ([2,t-1],t)^{(t-j)} \prec \mathrm{colex}((1, 2), (1, j))                                                         \\
   & \xrightarrow{\text{s.e.'s with $(1, j+1)$}} \bigprec_{3 \leq t \leq j} ([2,t-1],t)^{(0)} \prec ([2,j],j+1)^{(0)} \prec ([2,j],j+2)^{(0)}                                            \\
   & {\color{red}{\prec (1, j+1) \prec (1, j+2) \prec (j+1, j+2)}} \prec \bigprec_{j+3 \leq t \leq n} ([2,t-1],t)^{(t-j)} \prec \mathrm{colex}((1, 2), (1, j))                           \\
   & \xrightarrow{\text{flip}} \bigprec_{3 \leq t \leq j} ([2,t-1],t)^{(0)} \prec ([2,j],j+1)^{(0)} \prec ([2,j],j+2)^{(0)}                                                              \\
   & \prec (j+1, j+2) \prec (1, j+2) \prec (1, j+1) \prec \bigprec_{j+3 \leq t \leq n} ([2,t-1],t)^{(t-j)} \prec \mathrm{colex}((1, 2), (1, j))                                          \\
\end{align*}

(VI) Order changes when $x(P_{234}) - \epsilon < t < x(P_{1n-1n}) - \epsilon$.
\begin{align*}
   & {\color{red}{(2, 3) \prec (2, 4) \prec (3, 4)}} \prec \bigprec_{5 \leq t \leq n} ([2,t-1],t)^{(0)} \prec \mathrm{colex}((1, 2), (1, n))            \\
   & \xrightarrow{\text{flip}} (3, 4) \prec (2, 4) \prec (2, 3) \prec \bigprec_{5 \leq t \leq n} ([2,t-1],t)^{(0)} \prec \mathrm{colex}((1, 2), (1, n)) \\
\end{align*}

The skeleton associated with the wiring diagram can be uniquely reconstructed by tracing the sequence of order changes in reverse.
More precisely, we traverse the intersection points in the reverse of their creation order.
At each step, we perform the following local operation on the corresponding intersection edges:
\begin{itemize}
  \item \textbf{Double point $P_{ij|kl}$.}
        The two endpoints $l_{ij} \cap l_t$ and $l_{kl} \cap l_t$ are rotated clockwise by a half twist around the segment connecting them.
  \item \textbf{Triple point $P_{ijk}$.}
        The three points $l_{ij} \cap l_t$, $l_{ik} \cap l_t$, and $l_{jk} \cap l_t$ are rotated clockwise by a half twist around the segment joining them, with $l_{ik} \cap l_t$ as the center of rotation.
\end{itemize}
These local operations are uniquely determined by the local geometry near each intersection point.
Hence, tracing the intersections in the reverse of their creation order yields a uniquely determined skeleton.

Furthermore, when the intersection points $l_{ij} \cap l_t$ and $l_{kl} \cap l_t$ in the skeleton correspond to the relation $[g_{ij}, g_{kl}] = 1$, the associated geometric operation has the following interpretation.
Geometrically, applying this relation corresponds to rotating the segment joining $l_{ij} \cap l_t$ and $l_{kl} \cap l_t$ counterclockwise by a full twist.
Taking into account this correspondence between algebraic relations and geometric transformations, one arrives at the final form of the skeletons depending on the relative positions of the indices, three types of skeletons occur (see Figure~\ref{fig:skeletons} for the corresponding skeletons):
\begin{itemize}
  \item non--crossing type, when $1 \leq i < j < r < s \leq n$ or $1 \leq i < r < s < j \leq n$;
  \item crossing type, when $1 \leq i < r < j < s \leq n$;
  \item admissible type, when $1 \leq i < j < k \leq n$.
\end{itemize}

\begin{figure}[h]
  \centering
  \begin{minipage}{0.29\hsize}
    \centering
    \begin{tikzpicture}
      \fill (0,0.5) circle (2pt);
      \fill (0,2) circle (2pt);
      \node at (0,0) {$\vdots$};
      \node at (0.5,0.5) {$(r, s)$};
      \node at (0,1.3) {$\vdots$};
      \node at (0.5,2) {$(i, j)$};
      \node at (0,2.5) {$\vdots$};
      \draw (0,0.5) .. controls (-0.3,1.3) .. (0,2);
    \end{tikzpicture}
  \end{minipage}
  \begin{minipage}{0.29\hsize}
    \centering
    \begin{tikzpicture}[scale=1]
      \fill (0,0) circle (2pt);
      \fill (0,1.5) circle (2pt);
      \fill (0,3) circle (2pt);
      \node at (0,-0.3) {$\vdots$};
      \node at (0.5,0) {$(r, s)$};
      \node at (0,0.6) {$\vdots$};
      \node at (0.8,1.5) {$(i, s)$};
      \node at (0,2.6) {$\vdots$};
      \node at (0.5,3) {$(i, j)$};
      \node at (0,3.5) {$\vdots$};
      \draw
      (0,0)
      to[out=170, in=-110, looseness=0.8] (0.1,1.3)
      to[out=70, in=-160, looseness=0.8] (0,3);
    \end{tikzpicture}
  \end{minipage}
  \begin{minipage}{0.29\hsize}
    \centering
    \begin{tikzpicture}
      \fill (0,0.5) circle (2pt);
      \fill (0,1.5) circle (2pt);
      \fill (0,2.5) circle (2pt);
      \node at (0,0) {$\vdots$};
      \node at (0.5,0.5) {$(j, k)$};
      \node at (0,1.1) {$\vdots$};
      \node at (0.5,1.5) {$(i, k)$};
      \node at (0,2.1) {$\vdots$};
      \node at (0.5,2.5) {$(i, j)$};
      \node at (0,3.1) {$\vdots$};
      \draw (0,0.5) .. controls (-0.3,1) .. (0,1.5);
      \draw (0,1.5) .. controls (-0.3,2) .. (0,2.5);
    \end{tikzpicture}
  \end{minipage}
  \caption{Skeletons for the non--crossing, crossing, and admissible types, from left to right.}\label{fig:skeletons}
\end{figure}
Thus, we obtain the following theorem.
\begin{thm}\label{thm:PBn_pres}
  Let $PB_n = \pi_1(\Cb^n \setminus Br(n))$.
  Then, $PB_n$ has a presentation with generators $g_{ij}$, $1 \leq i < j \leq n$, and relations given by the following:
  \begin{itemize}
    \item[(r-I)] $[g_{ij}, g_{rs}] = 1$, where $1 \leq i < j < r < s \leq n$ or $1 \leq i < r < s < j \leq n$.
    \item[(r-II)] $[g_{ij}, g_{is} g_{rs} g_{is}^{-1}] = 1$, where $1 \leq i < r < j < s \leq n$.
    \item[(r-III)] $[g_{ij}, g_{ik}g_{jk}] = [g_{ij} g_{ik}, g_{jk}] = 1$, where $1 \leq i < j < k \leq n$,
  \end{itemize}
  where $g_{ij}$ denotes the meridian corresponding to the line $l_{ij}$.
\end{thm}
As a consequence of the above computation, we obtain the following corollary.
\begin{cor}\label{cor:PB_equals_MM}
  Two presentations of $PB_n$ and $PB_K$ are identical via $g_{ij} \longmapsto S_{ij}$.
\end{cor}
\section{Manin--Schechtman arrangements}\label{sec:MS_pi1}
In this section, we generalize the computation of the fundamental group of the complement of the braid arrangement to the Manin--Schechtman arrangements.

We begin by recalling the definition of the Manin--Schechtman arrangements \cite{MS89}.
Fix a central generic hyperplane arrangement $\Ac^0 = \{ H_1^0, H_2^0, \dots, H_n^0 \}$ in $\Cb^k$, such that any $m \leq k$ hyperplanes intersect in codimension $m$, and the only common intersection of all hyperplanes is the origin.
Consider the space of all parallel translates of $\Ac^0$,
\begin{equation*}
  \Sb = \Sb(H_1^0, H_2^0, \dots, H_n^0) = \{ (H_1^{t_1}, H_2^{t_2}, \dots, H_n^{t_n}) \mid t_1, t_2, \dots, t_n \in \Cb \} ,
\end{equation*}
where $H_i^{t_i}$ denotes the affine hyperplane parallel to $H_i^0$ obtained by translating $H_i^0$ in the normal direction by the parameter $t_i$.
Via the correspondence $(H_1^{t_1}, H_2^{t_2}, \dots, H_n^{t_n}) \mapsto (t_1, t_2, \dots, t_n)$, we identify $\Sb$ with $\Cb^n$.

For each subset $L \subset [n]$ with $|L| = k + 1$, define
\begin{equation*}
  D_L = \Bigl \{( H_1^{t_1}, H_2^{t_2}, \dots, H_n^{t_n}) \in \Sb\ \Big|\ \bigcap_{p \in L} H_p^{t_p} \neq \emptyset \Bigr \}.
\end{equation*}
Equivalently, $D_L$ consists of those translates for which the subarrangement indexed by $L$ fails to be in general position.
In particular, $D_L$ is a hyperplane in $\Sb \simeq \Cb^n$.
The \textbf{Manin--Schechtman arrangement} is the hyperplane arrangement
\begin{equation*}
  MS(n, k) = \{ D_L \}_{\substack{L \subset [n]\\ |L| = k+1}},
\end{equation*}
introduced in \cite{MS89} as a higher analogue of the braid arrangement.
In particular, $MS(n, 1)$ coincides with the braid arrangement $Br(n)$.

It is known that the combinatorial structure of the Manin--Schechtman arrangement may depend on the choice of the underlying generic arrangement $\Ac^0$.
Following the terminology introduced by Bayer and Brandt~\cite{BB97}, the arrangement $\Ac^0$ is called very generic if this combinatorics is independent of the particular choice of $\Ac^0$, and non--very generic otherwise.

Throughout the present paper, we restrict ourselves to the very generic case.
Accordingly, we simply write $MS(n, k)$ for the corresponding Manin--Schechtman arrangement, suppressing the dependence on $\Ac^0$.
For discussions of the non--very generic case, see for instance~\cite{BB97, Ath99, Fal94, Sai25, Sai26}.

We now consider the family of line arrangements obtained as generic $2$--sections of $MS(n, k)$.
Note that in the very generic case, the arrangement $MS(n, k)$ has only two types of codimension-2 intersections,
\begin{align*}
  D_I \cap D_J, \quad I, J \subset [n], \quad |I \cap J| < k, & \quad \text{and} \\
  \bigcap_{j = 1}^{k + 2} D_{L \setminus \{ i_j \}}, \quad L = \{ i_1 < i_2 < \dots < i_{k + 2} \} \subset [n], \quad |L| = k + 2.
\end{align*}
The former gives the double points and the latter gives the $(k + 2)$--multiple points in the section.
\begin{defi}
  An \textbf{$(n, k)$--lexicographic line arrangement} is a real line arrangement consisting of lines $l_I$ indexed by $(k + 1)$--subsets $I = \{ i_1 < i_2 < \dots < i_{k+1} \} \subset [n]$, satisfying the following properties.
  \begin{enumerate}
    \item[(L1)] Let $C$ be a circle containing all intersection points of the lines $l_I$ in its interior.
          For each $I$, let $C_I = C \cap l_I$ be the intersection point lying on the ray of $l_I$ in the direction of increasing $x$--coordinate.
          Then the points $C_I$ appear in lexicographic order with respect to the indices $I$ when going along $C$ clockwise starting from $C_{\{ 1, 2, \dots, k + 1 \}}$.
    \item[(L2)] The intersection points are of exactly the following two types:
          \begin{itemize}
            \item double point $P_{I \mid J} = l_I \cap l_J$ when $|I \cap J| < k$ (equivalently, when $I$ and $J$ do not belong to the same $(k+2)$--packet in the sense of~\cite{MS89});
            \item $(k + 2)$--multiple point $P_L = \bigcap_{j = 1}^{k + 2} l_{L \setminus \{ i_j \}}$ associated with a $(k+2)$--subset $L = \{ i_1 < i_2 < \dots < i_{k + 2} \} \subset [n]$.
          \end{itemize}
    \item[(L3)] The $(k + 2)$--multiple points are totally ordered by their $x$--coordinates so that $x(P_{L'}) < x(P_L)$ whenever $L < L'$ in lexicographic order.
  \end{enumerate}
\end{defi}
Note that an $(n, 1)$--lexicographic line arrangement coincides with the lexicographic line arrangement of the braid arrangement introduced in Subsection~\ref{subsec:lex}.
Although the existence of an $(n, k)$--lexicographic line arrangement is implicit in the original work of Manin and Schechtman~\cite{MS89}, no proof is given there.
We therefore provide below a proof of existence before proceeding to the braid monodromy computation.
\begin{prop}\label{prop:nk-line-exist}
  For any positive integers $n, k$ with $n \geq k + 1$, there exists an $(n, k)$--lexicographic line arrangement.
\end{prop}
\begin{proof}
  Since the statement is trivial for $n = k + 1$, we assume $n \geq k + 2$.

  Choose a real number $R>2$, and choose pairwise distinct real numbers
  \begin{equation*}
    \lambda_1 > \lambda_2 > \cdots > \lambda_n > 0,
  \end{equation*}
  such that $\lambda_i \geq R \lambda_{i+1}$ for all $i$.
  This separation condition implies that the sums $\sum_{i \in I} \lambda_i$ are totally ordered by the lexicographic order on $(k + 1)$--subsets $I \subset [n]$.

  For $(x, y) \in \Rb^2$ set
  \begin{equation*}
    t_i(x, y) := \lambda_i^{\,k} x + \lambda_i^{\,k+1}y + c_i
    \qquad(1 \leq i \leq n).
  \end{equation*}
  For each $(k + 1)$--subset $I = \{ i_0 <i_1 < \cdots <i_k \} \subset [n]$, define the line $l_I \subset \Rb^2$ by
  \begin{equation*}
    \det
    \begin{pmatrix}
      1      & \lambda_{i_0} & \lambda_{i_0}^2 & \cdots & \lambda_{i_0}^{k-1} & t_{i_0}(x, y) \\
      1      & \lambda_{i_1} & \lambda_{i_1}^2 & \cdots & \lambda_{i_1}^{k-1} & t_{i_1}(x, y) \\
      \vdots & \vdots        & \vdots          &        & \vdots              & \vdots        \\
      1      & \lambda_{i_k} & \lambda_{i_k}^2 & \cdots & \lambda_{i_k}^{k-1} & t_{i_k}(x, y)
    \end{pmatrix}=0.
  \end{equation*}
  Expanding along the last column, we obtain an equation
  \begin{equation*}
    A_I x + B_I y + D_I = 0,
  \end{equation*}
  where $D_I$ depends linearly on $(c_1,\dots,c_n)$ and
  \begin{equation*}
    A_I = V_I, \quad B_I = V_I \sum_{i \in I} \lambda_i,
  \end{equation*}
  with $V_I \neq 0$ the Vandermonde factor.
  Hence the slope is
  \begin{equation*}
    m_I = - \frac{1}{\sum_{i \in I} \lambda_i}.
  \end{equation*}
  Since the strong separation condition on the $\lambda_i$ implies that
  $I < J$ in lexicographic order yields
  $\sum_{i \in I} \lambda_i > \sum_{j \in J} \lambda_j$,
  the slopes $m_I$ are ordered lexicographically.
  Thus, condition~(L1) follows.

  Next, a point $(x, y)$ lies on $l_I$ if and only if there exists a polynomial $p(X)$ of degree at most $k - 1$ such that $t_i(x, y) = p(\lambda_i)$ for all $i \in I$.
  For a $(k + 2)$--subset $L$, this condition for all $I = L \setminus \{ i \}$ is equivalent to a square linear system whose coefficient matrix is of Vandermonde type.
  Since the $\lambda_i$ are distinct, this system has a unique solution, defining a unique $(k + 2)$--multiple point $P_L$.

  Any unwanted multiple intersection or coincidence $P_L = P_{L'}$ is described by explicit rank conditions, hence by polynomial equations in the parameters $c_i$.
  Therefore the set $U_{(L2)} \subset \Rb^n$ of parameters for which only the prescribed double points and $(k + 2)$--multiple points occur is a nonempty Zariski open subset.

  Finally, the coordinates of each $P_L$ depend affinely on the parameters $c_i$.
  Ordering the $(k + 2)$--subsets lexicographically and comparing consecutive $x$--coordinates yields finitely many linear inequalities.
  Since the corresponding linear functions are not identically zero, these inequalities define a nonempty open subset $U_{(L3)} \subset \Rb^n$.
  Choosing $(c_1, c_2, \dots, c_n) \in U_{(L2)} \cap U_{(L3)}$ yields an arrangement satisfying (L2) and (L3).

  Together with (L1), this completes the proof.
\end{proof}

Below, we present the computation that yields a presentation of the Manin--Schechtman group.
The lengthy computation in Subsection~\ref{subsec:monodromy_tech} describes in detail the evolution of the vertical order of sectional intersection points for the braid arrangement $Br(n) = MS(n, 1)$.
From a combinatorial viewpoint, this calculation serves as a prototype for the higher Manin--Schechtman arrangements $MS(n, k)$.
Indeed, the braid monodromy algorithm depends only on the combinatorics of the $(k + 2)$--multiple points and the double points in the induced wiring diagram, and therefore extends formally to $MS(n, k)$.
For simplicity, we treat explicitly only the case $k = 2$; however, by the above reasoning, the same method applies verbatim to the general case $k > 2$.

Conceptually, passing from $MS(n, 1)$ to $MS(n, 2)$ amounts to replacing pairs $ij$ labeling the lines by triples $ijk$, and replacing each triple point $P_{ijk}$ by a $4$--multiple point corresponding to a $4$--subset $\{ i, j, k, l \}$.
The elementary local operations governing the order changes remain unchanged: single exchanges arising from double points, and flips corresponding to the multiple points.
Thus, one may apply the same Moishezon--Teicher procedure to a generic $(n, 2)$--lexicographic section of $MS(n, 2)$ to obtain a well--defined sequence of skeletons.
For brevity, we omit the explicit order--change computation in this case.
\begin{rem}
  The orderings constructed here provide a geometric realization of the higher Bruhat order $B(n, 2)$, in terms of real line arrangements and their intersection points.
\end{rem}

From the point of view of the fundamental group, only the resulting local skeletons are relevant.
Double points in the $(n, 2)$--lexicographic line arrangement produce commutativity relations analogous to (R-I), (R-II) and (R-III) for the pure braid group, now written in terms of the ``$3$--swings type'' and the ``$4$--swings type''.
Collecting all relations obtained from the skeletons gives a presentation of the ``$2$-dimensional pure braid group'' $\pi_1(\Cb^n \setminus MS(n, 2))$, generalising the modified Artin presentation.
As a summary, we have the following theorem.
\begin{thm}
  The fundamental group $\pi_1(\Cb^n \setminus MS(n, 2))$ of the complement of the Manin--Schechtman arrangement $MS(n,2)$ admits the presentation
  \begin{equation*}
    \pi_1(\Cb^n \setminus MS(n, 2)) = \langle\, g_{i_1i_2i_3} \mid \text{relations of types {\rm (R-I)}--{\rm (R-IV)}} \,\rangle,
  \end{equation*}
  where the relations are described as follows.
  Assume that $1 \leq i_1 < i_2 < i_3 \leq n$, $1 \leq j_1 < j_2 < j_3 \leq n$, and $i_1 \leq j_1$.
  \begin{itemize}
    \item[(R-I)]
          The generators $g_{i_1i_2i_3}$ and $g_{j_1j_2j_3}$ commute,
          \begin{equation*}
            [g_{i_1i_2i_3}, g_{j_1j_2j_3}] = 1,
          \end{equation*}
          whenever the pair of triples $(i_1,i_2,i_3)$ and $(j_1,j_2,j_3)$ does not fall into any of the cases (R-II) or (R-III).
    \item[(R-II)]
          Suppose that the two triples share exactly two indices.
          According to the position of the common indices, one obtains the following three types of relations:
          \begin{itemize}
            \item $[g_{i_1i_2i_3}, g_{i_1i_2j_3} g_{i_1j_2j_3} g_{i_1i_2j_3}^{-1}] = 1$ when $1 \leq i_1 = j_1 < i_2 < j_2 < i_3 < j_3 \leq n$;
            \item $[g_{i_1i_2i_3}, g_{i_1i_2j_3} g_{j_1i_2j_3} g_{i_1i_2j_3}^{-1}] = 1$ when $1 \leq i_1 < j_1 < i_2 = j_2 < i_3 < j_3 \leq n$;
            \item $[g_{i_1i_2i_3}, g_{i_1j_2i_3} g_{j_1j_2i_3} g_{i_1j_2i_3}^{-1}] = 1$ when $1 \leq i_1 < j_1 < i_2 < j_2 < j_3 = i_3 \leq n$.
          \end{itemize}
          Figures~\ref{fig:rel2-1}, \ref{fig:rel2-2} and \ref{fig:rel2-3} illustrate the corresponding skeletons and configuration of indices.
    \item[(R-III)] For $1 \leq i_1 < j_1 \leq i_2 < j_2 \leq i_3 < j_3 \leq n$, the following relation holds:
          \begin{equation*}
            [g_{i_1i_2i_3}, g_{i_1j_2j_3} g_{i_1i_2j_3} g_{j_1j_2j_3} g_{i_1i_2j_3}^{-1} g_{i_1j_2j_3}^{-1}] = 1.
          \end{equation*}
          Figure~\ref{fig:rel3} illustrates the corresponding skeleton and configuration of indices.
    \item[(R-IV)]For $1 \leq i_1 < i_2 < i_3 < i_4 \leq n$, the following three relations hold:
          \begin{align*}
             & [g_{i_1i_2i_3}, g_{i_1i_2i_4} g_{i_1i_3i_4} g_{i_2i_3i_4}] = [g_{i_1i_2i_3} g_{i_1i_2i_4}, g_{i_1i_3i_4} g_{i_2i_3i_4}] \\
             & = [g_{i_1i_2i_3} g_{i_1i_2i_4} g_{i_1i_3i_4}, g_{i_2i_3i_4}] = 1.
          \end{align*}
          Figure~\ref{fig:rel4} shows the associated skeleton and configuration of indices.
  \end{itemize}
\end{thm}

\begin{figure}[h]
  \centering
  \begin{minipage}{0.48\linewidth}
    \centering
    \begin{tabular}{cc}
      \begin{minipage}{0.35\linewidth}
        \centering
        \begin{tikzpicture}[scale=0.8]
          \fill (0,0.5) circle (2pt);
          \fill (0,1.5) circle (2pt);
          \fill (0,2.5) circle (2pt);
          \node at (0,0) {$\vdots$};
          \node at (1.3,0.5) {$(j_1, j_2, j_3)$};
          \node at (0,2.1) {$\vdots$};
          \node at (1.3,1.5) {$(i_1, i_2, j_3)$};
          \node at (0,3.1) {$\vdots$};
          \node at (1.3,2.5) {$(i_1, i_2, i_3)$};
          \node at (0,1.1) {$\vdots$};
          \draw (0,0.5) .. controls (-0.5,1.3) and (0.2,1.3) .. (0.2,1.5);
          \draw (0.2,1.5) .. controls (0.2,1.7) and (-0.5,1.7) .. (0,2.5);
        \end{tikzpicture}
      \end{minipage}
       &
      \begin{minipage}{0.45\linewidth}
        \centering
        \begin{tikzpicture}[scale=0.8]
          \begin{scope}[scale=0.8]
            \node at ({-sqrt(3)-1},1) {$i_1 = j_1$};
            \node at (0,-2.3) {$i_2$};
            \node at ({sqrt(3)+0.3},-1) {$j_2$};
            \node at ({sqrt(3)+0.3},1) {$i_3$};
            \node at (0,2.3) {$j_3$};
            \draw ({-sqrt(3)},1) -- (0,-2) -- ({sqrt(3)},1) -- cycle;
            \draw ({-sqrt(3)},1) -- ({sqrt(3)},-1) -- (0,2) -- cycle;
          \end{scope}
        \end{tikzpicture}
      \end{minipage}
    \end{tabular}
    \medskip
    \captionof{figure}{Skeleton and configuration of indices yielding a first relation of type (R-II). \label{fig:rel2-1}}
  \end{minipage}
  \hfill
  \begin{minipage}{0.48\linewidth}
    \centering
    \begin{tabular}{cc}
      \begin{minipage}{0.35\linewidth}
        \centering
        \begin{tikzpicture}[scale=0.8]
          \fill (0,0.5) circle (2pt);
          \fill (0,1.5) circle (2pt);
          \fill (0,2.5) circle (2pt);
          \node at (0,0) {$\vdots$};
          \node at (1.3,0.5) {$(j_1, j_2, j_3)$};
          \node at (0,2.1) {$\vdots$};
          \node at (1.3,1.5) {$(i_1, j_2, i_3)$};
          \node at (0,3.1) {$\vdots$};
          \node at (1.3,2.5) {$(i_1, i_2, i_3)$};
          \node at (0,1.1) {$\vdots$};
          \draw (0,0.5) .. controls (-0.5,1.3) and (0.2,1.3) .. (0.2,1.5);
          \draw (0.2,1.5) .. controls (0.2,1.7) and (-0.5,1.7) .. (0,2.5);
        \end{tikzpicture}
      \end{minipage}
       &
      \begin{minipage}{0.45\linewidth}
        \centering
        \begin{tikzpicture}[scale=0.8]
          \begin{scope}[scale=0.8]
            \node at ({-sqrt(3)-0.3},1) {$i_1$};
            \node at ({-sqrt(3)-0.3},-1) {$j_1$};
            \node at (0,-2.3) {$i_2 = j_2$};
            \node at ({sqrt(3)+0.3},1) {$i_3$};
            \node at (0,2.3) {$j_3$};
            \draw ({-sqrt(3)},1) -- (0,-2) -- ({sqrt(3)},1) -- cycle;
            \draw ({-sqrt(3)},-1) -- (0,-2) -- (0,2) -- cycle;
          \end{scope}
        \end{tikzpicture}
      \end{minipage}
    \end{tabular}
    \medskip
    \captionof{figure}{Skeleton and configuration of indices yielding a second relation of type (R-II). \label{fig:rel2-2}}
  \end{minipage}
\end{figure}

\begin{figure}[h]
  \centering
  \begin{minipage}{0.48\linewidth}
    \centering
    \begin{tabular}{cc}
      \begin{minipage}{0.35\linewidth}
        \centering
        \begin{tikzpicture}[scale=0.8]
          \fill (0,0.5) circle (2pt);
          \fill (0,1.5) circle (2pt);
          \fill (0,2.5) circle (2pt);
          \node at (0,0) {$\vdots$};
          \node at (1.3,0.5) {$(j_1, j_2, j_3)$};
          \node at (0,2.1) {$\vdots$};
          \node at (1.3,1.5) {$(i_1, j_2, i_3)$};
          \node at (0,3.1) {$\vdots$};
          \node at (1.3,2.5) {$(i_1, i_2, i_3)$};
          \node at (0,1.1) {$\vdots$};
          \draw (0,0.5) .. controls (-0.5,1.3) and (0.2,1.3) .. (0.2,1.5);
          \draw (0.2,1.5) .. controls (0.2,1.7) and (-0.5,1.7) .. (0,2.5);
        \end{tikzpicture}
      \end{minipage}
       &
      \begin{minipage}{0.45\linewidth}
        \centering
        \begin{tikzpicture}[scale=0.8]
          \begin{scope}[scale=0.8]
            \node at ({-sqrt(3)-0.3},1) {$i_1$};
            \node at ({-sqrt(3)-0.3},-1) {$j_1$};
            \node at (0,-2.3) {$i_2$};
            \node at ({sqrt(3)+0.3},-1) {$j_2$};
            \node at ({sqrt(3)},1.3) {$i_3=j_3$};
            \draw ({-sqrt(3)},1) -- (0,-2) -- ({sqrt(3)},1) -- cycle;
            \draw ({-sqrt(3)},-1) -- ({sqrt(3)},-1) -- ({sqrt(3)},1) -- cycle;
          \end{scope}
        \end{tikzpicture}
      \end{minipage}
    \end{tabular}
    \medskip
    \captionof{figure}{Skeleton and configuration of indices yielding a third relation of type (R-II).\label{fig:rel2-3}}
  \end{minipage}
  \hfill
  \begin{minipage}{0.48\linewidth}
    \centering
    \begin{tabular}{cc}
      \begin{minipage}{0.35\linewidth}
        \centering
        \begin{tikzpicture}[scale=0.8]
          \fill (0,0.5) circle (2pt);
          \fill (0,1.5) circle (2pt);
          \fill (0,2.5) circle (2pt);
          \node at (0,0) {$\vdots$};
          \node at (1.3,0.5) {$(j_1, j_2, j_3)$};
          \node at (0,2.1) {$\vdots$};
          \node at (1.3,1.5) {$(i_1, i_2, j_3)$};
          \node at (0,3.1) {$\vdots$};
          \node at (1.3,2.5) {$(i_1, i_2, i_3)$};
          \node at (0,1.1) {$\vdots$};
          \draw (0,0.5) .. controls (-0.5,1.3) and (0.2,1.3) .. (0.2,1.5);
          \draw (0.2,1.5) .. controls (0.2,1.7) and (-0.5,1.7) .. (0,2.5);
        \end{tikzpicture}
      \end{minipage}
       &
      \begin{minipage}{0.45\linewidth}
        \centering
        \begin{tikzpicture}[scale=0.8]
          \begin{scope}[scale=0.8]
            \node at ({-sqrt(3)-0.3},1) {$i_1$};
            \node at ({-sqrt(3)-0.3},-1) {$j_1$};
            \node at (0,-2.3) {$i_2$};
            \node at ({sqrt(3)+0.3},-1) {$j_2$};
            \node at ({sqrt(3)+0.3},1) {$i_3$};
            \node at (0,2.3) {$j_3$};
            \draw ({-sqrt(3)},1) -- (0,-2) -- ({sqrt(3)},1) -- cycle;
            \draw ({-sqrt(3)},-1) -- ({sqrt(3)},-1) -- (0,2) -- cycle;
          \end{scope}
        \end{tikzpicture}
      \end{minipage}
    \end{tabular}
    \medskip
    \captionof{figure}{Skeleton and configuration of indices yielding a relation of type (R-III).\label{fig:rel3}}
  \end{minipage}
\end{figure}

\begin{figure}[h]
  \centering
  \begin{minipage}[b]{0.2\linewidth}
    \centering
    \begin{tikzpicture}[scale=0.8]
      \fill (0,0.5) circle (2pt);
      \fill (0,1.5) circle (2pt);
      \fill (0,2.5) circle (2pt);
      \fill (0,3.5) circle (2pt);
      \node at (0,0) {$\vdots$};
      \node at (1.3,0.5) {$(j_1, j_2, j_3)$};
      \node at (0,1.1) {$\vdots$};
      \node at (1.3,1.5) {$(i_1, i_2, i_4)$};
      \node at (0,2.1) {$\vdots$};
      \node at (1.3,2.5) {$(i_1, i_3, i_4)$};
      \node at (0,3.1) {$\vdots$};
      \node at (1.3,3.5) {$(i_2, i_3, i_4)$};
      \node at (0,4.2) {$\vdots$};
      \draw (0,0.5) .. controls (-0.5,1) .. (0,1.5);
      \draw (0,1.5) .. controls (-0.5,2) .. (0,2.5);
      \draw (0,2.5) .. controls (-0.5,3) .. (0,3.5);
    \end{tikzpicture}
  \end{minipage}
  \begin{minipage}[b]{0.2\linewidth}
    \centering
    \begin{tikzpicture}[scale=0.8]
      \node at (0,-0.3) {$i_1$};
      \node at (2,0) {$i_2$};
      \node at (1.7,2.3) {$i_3$};
      \node at (0,2.1) {$i_4$};
      \draw (0,0) -- (2,0.3) -- (1.7,2) -- (0,1.8) -- cycle;
    \end{tikzpicture}
  \end{minipage}
  \captionof{figure}{Skeleton and configuration of indices yielding a relation of type (R-IV).\label{fig:rel4}}
\end{figure}
From an enumerative point of view, the above presentation involves $\binom{n}{3}$ generators with $\binom{\binom{n}{3}}{2} - \binom{n}{6} - 5 \binom{n}{5} - \binom{n}{4}$ numbers of relations of type (R-I), $3 \binom{n}{5}$ numbers of relations of type (R-II), $\binom{n}{6} + 2 \binom{n}{5} + \binom{n}{4}$ numbers of relations of type (R-III), and $3 \binom{n}{4}$ numbers of relations of type (R-IV).

\begin{rem}
  Lawrence~\cite{Law91} studied the Manin--Schechtman (higher braid) groups and noted that, for $k > 1$, these groups admit no unique natural presentation.
  The analysis in~\cite{Law91} yields presentations depending on auxiliary choices, reflecting this inherent non-uniqueness.
  From this viewpoint, it would be interesting to compare Lawrence's approach with the presentation obtained in the present paper.

  Related notions of higher braid groups, such as the groups $G_n^k$ introduced
  by Manturov and Nikonov~\cite{MN15}, may also provide an interesting point of
  comparison.
\end{rem}
\vskip\baselineskip
\textit{Acknowledgements.}
The author is grateful to Anatoly Libgober for drawing his attention to the work~\cite{MM09}, and for the opportunity to visit the University of Illinois at Chicago, where part of this work was developed.
The author also thanks Toshiyuki Akita for carefully reading an earlier draft of the manuscript and for providing helpful comments.
This work was supported by JSPS KAKENHI Grant Number JP24K16926 and JPJSBP120256504.
\bibliography{pre_pbn}
\bibliographystyle{amsalpha}
\end{document}